\documentclass{amsart}
\usepackage{amssymb, amscd}
\usepackage{pst-all}
\numberwithin{equation}{section}

\def\FF{{\mathbb F}}

\def\ZZ{{\mathbb Z}}

\newcommand\proofsquare{\nobreak\hfill \hbox{%
\vrule height 5pt
\kern-.4pt
 \vbox{%
\hrule width 5pt depth0pt height.4pt
 \kern4.6pt \hrule  }
\kern-3.75pt
\vrule height 5pt}\kern1pt
\par}

\newtheorem{definition-lemma}[theorem]{Definition-Lemma}

\theoremstyle{definition}

\theoremstyle{remark}

\begin{document}
\title[Hunting for Curves with Many Points]{Hunting for Curves with Many Points}
\author{Gerard van der Geer}

\address{Korteweg-de Vries Instituut, Universiteit van Amsterdam, Plantage Muidergracht 24, 1018 TV Amsterdam, The Netherlands}

\email{geer@science.uva.nl}

\subjclass{14K10}

\begin{abstract}
We construct curves with many points over finite fields by using
the class group.
\end{abstract}
\maketitle

\begin{section}{Introduction}
\label{sec: intro}
The question how many rational points a curve of given genus
over a finite field of given cardinality can have is an attractive challenge
and the table of curves with many points (\cite{Tables})
has been a reference point for progress for genus $\leq 50$ and small
fields of characteristic $p=2$ and $p=3$. 
The tables record there for a pair $(q,g)$ an interval $[a,b]$ where
$b$ is the best upper bound for the maximum number of points of a 
curve of genus $g$ over ${\FF}_q$ and $a$ gives a lower bound
obtained from an explicit example of a curve $C$ defined over ${\FF}_q$ with 
$a$ (or at least $a$) rational points. So $N_q(g) \in [a,b]$ with
$N_q(g)$ the maximum number of rational points on a smooth connected 
projective curve defined over ${\FF}_q$.
For progress on the upper bounds we
refer to the papers of Howe and Lauter, \cite{HL}  
and to the references given in the tables.
It is the purpose of this little
paper to record recent improvements of the table and at the same time
give constructions for many of the present records of the table.
The methods employ the class group and are variations on 
well-known themes and do not involve new ideas.
We hope it inspires people to take up the challenge to improve the tables.
\end{section}
\begin{section}{Using the Class Group}
Let $C$ be a curve of genus $g$ defined over $\FF_2$ with `many' rational points.
If $f$ is a rational function not of the form $h^2+h$ for $h$ in 
the function field ${\FF}_2(C)$ 
of $C$ then in the Artin-Schreier cover
of $C$ defined by $w^2+w=f$ all the rational points of $C$ where $f$ vanishes will
be split. If the contribution of the poles of $f$ is limited and
$f$ vanishes in many points this may yield
curves over $\FF_2$ with many points. Serre applied this method 
(cf.\ \cite{S1})
to construct curves over ${\FF}_2$ of genus $g=11$ with $14$ rational points
(resp.\ $g=13$ with $N=15$ and $g=10$ with $N=13$) thus determining $N_2(10)=13$,
$N_2(11)=14$ and $N_2(13)=15$. To find suitable $f$ one uses the class group.
Serre used an elliptic curve, but we may of course use higher genus curves as well.
We illustrate the 
method by constructing many curves that reach or improve the present
`records'.

\begin{subsection}{Base Curve of Genus $1$ over ${\FF}_2$}
Let $C$ be the (projective smooth) curve of genus $1$ defined 
by the affine equation $y^2+y=x^3+x$. It has
five rational points $P_0=\infty$, $P_1=(0,0)$, $P_2=(1,0)$, $P_3=(1,1)$ and $P_4=(0,1)$.
Since we have an isomorphism
$C(\FF_2)\cong {\ZZ}/5{\ZZ}$ of abelian groups such that under the isomorphism
$P_i$ corresponds to $i \, (\bmod \, 5)$, the divisor $\sum a_iP_i$ is of degree $0$ and linearly equivalent to $0$ if and only if 
$$
\sum_{i=0}^4 a_i=0 \quad {\rm and} \quad \sum_{i=0}^4 a_i \, i \equiv
0\, (\bmod \, 5).\eqno(1)
$$ 
So there exists a function $f \in {\FF}_2(C)$ 
with divisor $(f)$ equal to a given divisor $\sum a_iP_i$
if and only if  the $5$-tuple 
$(a_0,a_1,\ldots,a_4)\in {\ZZ}^5$ satisfies (1).

Often we shall
denote a relation for linear equivalence 
$\sum_{i=0}^4 a_iP_i \sim 0$ simply by 
$[a_0,a_1,\ldots,a_4]$ for shortness. 
Two examples of such relations for the present curve
are $[-3,-1,2,1,1]$ and $[-1,-3,1,1,2]$ with
functions $f_1$ and $f_2$. This gives curves $C_{f_1}$ and $C_{f_2}$ 
of genera $g=4$ with $N=8$ rational points. Note that $N_2(4)=8$, 
so these curves realize the maximum for genus~$4$. 
The fibre product $C_{f_1}\times_{C}
C_{f_2}$ has genus $g=11$ with $N=14$ rational points. Since the upper bound
for $N_2(11)$ is $14$, Serre (cf.\ \cite{S1}) thus showed that $N_2(11)=14$. 

Another example is the relation 
$-7P_0+2P_1+3P_2+P_4+P_5 \sim 0$. If $f_3$ is a
corresponding function then the Artin-Schreier cover 
is a curve $C_{f_3}$ with
genus $g=5$ and with $N=9$ rational points as one readily checks. Note that
$N_2(5)=9$. The fibre product $C_{f_1}\times_C C_{f_3}$ is a curve of genus $g=13$
with $N=15$ rational points, hence $N_2(13)=15$, again due to Serre.

Combining $f_1$ with the function $f_4$ corresponding to 
the relation $[-3,2,1,-1,1]$
gives $C_{f_1}\times_CC_{f_4}$ of genus $10$ with $N=13$ rational points.
Note that the divisor of $f_1+f_4$ is 
$-2P_0-P_1+P_2+2P_4+P_5$. 
Artin-Schreier reduction shows that the conductor of this cover is 
$2P_0+2P_1+2P_3$. This determines $N_2(10)=13$, a result due to Serre.

We list here eight relations that can be used to obtain more good curves;
\smallskip

\vbox{
\bigskip\centerline{\def\quad{\hskip 0.6em\relax}
\def\quod{\hskip 0.5em\relax }
\vbox{\offinterlineskip
\hrule
\halign{&\vrule#&\strut\quod\hfil#\quad\cr
height2pt&\omit&&\omit&&\omit&&\omit&\cr
& n && relation  && n && relation &\cr
\noalign{\hrule}
& $1$ && $[-3,-1,2,1,1]$  && $5$ && $[-5,-1,3,2,1]$  & \cr
& $2$ && $[-1,-3,1,1,2]$  && $6$ && $[-9,1,3,2,3]$  & \cr
& $3$ && $[-7,2,3,1,1]$  && $7$ && $[-11,1,4,3,3]$  & \cr
& $4$ && $[-3,2,1,-1,1]$  && $8$ && $[-13,2,4,3,4]$  & \cr
} \hrule}
}}

\noindent
and reap the fruits by associating to a tuple $f_{i_1},\ldots,f_{i_r}$
the corresponding fibre product 
$C_{f_{i_1}}\times_C C_{f_2} \times \cdots \times_C C_{f_{i_r}}$.

\smallskip
\vbox{
\bigskip\centerline{\def\quad{\hskip 0.6em\relax}
\def\quod{\hskip 0.5em\relax }
\vbox{\offinterlineskip
\hrule
\halign{&\vrule#&\strut\quod\hfil#\quad\cr
height2pt&\omit&&\omit&&\omit&&\omit&&\omit&&\omit&&\omit&&\omit&\cr
& space && $g$ && $\# C({\FF}_2$ && interval &&  space && $g$ && $\# $ && interval &\cr
\noalign{\hrule}
& $\langle f_1\rangle$ && $4$ && $8$ && $[8]$ && $\langle f_1, f_3, f_5\rangle$ && $29$ && $25$ && $[25-27]$ &\cr
& $\langle f_3\rangle$ && $5$ && $9$ && $[9]$ && $\langle f_1, f_3, f_4\rangle$ && $30$ && $25$ && $[25-27]$ &\cr
& $\langle f_1, f_4\rangle$ && $10$ && $13$ && $[13]$ && $\langle f_1, f_2, f_3\rangle$ && $32$ && $27$ && $[26-29]$ &\cr
& $\langle f_1, f_2\rangle$ && $11$ && $14$ && $[14]$ && $\langle f_2, f_3, f_5\rangle$ && $34$ && $27$ && $[27-30]$ &\cr
& $\langle f_1, f_3\rangle $ && $13$ && $15$ && $[15]$ &&$\langle f_1, f_3, f_6\rangle$ && $35$ && $29$ && $[29-31]$ &\cr
& $\langle f_3, f_5\rangle $ && $14$ && $15$ && $[15-16]$ && $\langle f_3, f_6, f_7\rangle $ && $39$ && $33$ && $[33]$ &\cr
& $\langle f_3, f_6\rangle$ && $15$ && $17$ && $[17]$ &&$\langle f_3, f_6, f_8\rangle $ && $43$ && $33$ && $[33-36]$ & \cr
& $\langle f_1, f_2, f_5\rangle $ && $28$ && $25$ && $[25-26]$ &&
 $\langle f_3, f_7, f_8\rangle $ && $45$ && $33$ && $[33-37]$ &\cr
} \hrule}
}}

\noindent
This reproduces many records of the tables and 
gives one improvement of the tables for curves over 
${\FF}_2$, namely for $g=32$. 
\end{subsection}
\begin{subsection}{Base of Genus $2$ over ${\FF}_2$}
Let $C$ be the curve of genus $2$ defined by the affine equation $y^2+y=(x^2+x)/(x^3+x^2+1)$
over $\FF_2$. It has six rational points $P_i$ ($i=1,\ldots,6$) and 
class group ${\rm Jac}(C)(\FF_2)
\cong \ZZ/19\ZZ$. Let $\phi: C \to {\rm Jac}(C)$ be the Abel-Jacobi map given by $Q \mapsto
Q -\deg(Q) \, P_1$. For a suitable numbering of the $P_i$ the images of the $P_i$ in 
${\rm Jac}(C)(\FF_2)\cong \ZZ/19\ZZ$ are as follows

\begin{footnotesize}
\smallskip
\vbox{
\bigskip\centerline{\def\quad{\hskip 0.6em\relax}
\def\quod{\hskip 0.5em\relax }
\vbox{\offinterlineskip
\hrule
\halign{&\vrule#&\strut\quod\hfil#\quad\cr
height2pt&\omit&&\omit&&\omit&&\omit&&\omit&&\omit&&\omit&\cr
&$i$ && $1$ && $2$ && $3$ && $4$ && $5$ && $6$ &\cr
\noalign{\hrule}
&$\phi(P_i)$ && $0$ && $1$ && $14$ && $6$ && $16$ && $4$& \cr
} \hrule}
}}
\end{footnotesize}

\noindent
From this table we see that we have the following linear equivalence 
$$
-3P_1-3P_2+2P_3+2P_4+P_5+P_6\sim 0. \eqno(2)
$$ If $f$ is the function with the left hand side of (2) as divisor the curve $C_f$ obtained as
the double cover of $C$ defined by $w^2+w=f$ has genus $g=7$ and $10$ rational points. This is
optimal. Similarly the relation
$$
-9P_1+2P_2+2P_3+P_4+2P_5+2P_6 \sim 0  
$$
gives rise to a curve with $g=8$ and $11$ rational points, again an optimal curve over $\FF_2$.
Let $Q$ be the divisor of zeros of the polynomial $x^3+x^2+1$ on $C$. 
Then $\phi(Q)=3$ and $-Q+\sum_{i=1}^6 P_i\sim 0$ which gives us a curve of 
genus $9$ with $12$ points, again an optimal curve. 

We list a number of divisors of functions $f_n$ ($n=1,\ldots,10$):

\vbox{
\bigskip\centerline{\def\quad{\hskip 0.6em\relax}
\def\quod{\hskip 0.5em\relax }
\vbox{\offinterlineskip
\hrule
\halign{&\vrule#&\strut\quod\hfil#\quad\cr
height2pt&\omit&&\omit&&\omit&&\omit&\cr
& n && relation  && n && relation &\cr
\noalign{\hrule}
& $1$ && $[-3,-3,2,2,1,1]$  && $6$ && $[-11,3,2,2,3,1]$  & \cr
& $2$ && $[-9,2,2,1,2,2]$  &&  $7$ && $[1,0,-5,2,1,1]$  & \cr
& $3$ && $[1,1,1,1,1,1]\sim Q$  && $8$ && $[3,-3,-3,1,1,1]$  & \cr
& $4$ && $[-1,-5,1,2,2,1]$  && $9$ && $[-7,0,2,2,2,1]$  &  \cr
& $5$ && $[-5,-1,2,1,1,2]$  && $10$ && $[-13,3,1,3,3,3]$  & \cr
} \hrule}
}}

\noindent
and we find as results:

\smallskip
\vbox{
\bigskip\centerline{\def\quad{\hskip 0.6em\relax}
\def\quod{\hskip 0.5em\relax }
\vbox{\offinterlineskip
\hrule
\halign{&\vrule#&\strut\quod\hfil#\quad\cr
height2pt&\omit&&\omit&&\omit&&\omit&&\omit&&\omit&&\omit&&\omit&\cr
& space && $g$ && $\# C({\FF}_2)$ && interval && space && $g$ && $\# C({\FF}_2)$ && interval &\cr
\noalign{\hrule}
& $\langle f_1\rangle $ && $7$ && $10$ && $[10]$ && $\langle f_2, f_5\rangle $ && $20$ && $19$ && $[19-21]$  &\cr
& $\langle f_2\rangle $ && $8$ && $11$ && $[11]$  && $\langle f_2, f_6\rangle $ && $22$ && $21$ && $[21-22]$  &\cr
& $\langle f_3\rangle $ && $9$ && $12$ && $[12]$  && $\langle f_2, f_{10} \rangle $ && $24$ && $21$ && $[21-23]$  &\cr
& $\langle f_4, f_8\rangle$ && $17$ && $17$ && $[17-18]$  && $\langle f_1, f_4, f_5\rangle $ && $43$ && $34$ && $[33-36]$ &\cr
& $\langle f_1, f_4\rangle$ && $18$ && $18$ && $[18-19]$  &&$\langle f_2, f_5, f_9\rangle $ && $44$ && $33$ && $[33-37]$ & \cr
} \hrule}
}}
Again, we find one improvement, namely for $g=43$.

We can also employ the class group for making an unramified
cover of $C$ of degree $19$ in which $P_1$ splits completely. 
The genus of this cover is $g=20$ and it contains $19$ rational points.
The interval is $[19-21]$. 
\end{subsection}
\begin{subsection}{Base of Genus $3$ over ${\FF}_2$}
Consider now the curve $C$ of genus $3$ defined over $\FF_2$ by the 
homogeneous equation
$$
x^3y+x^3z+x^2y^2+xz^3+y^3z+y^2z^2=0.
$$
It has $7$ rational points and its class group ${\rm Jac}(C)(\FF_2)$
is isomorphic to $\ZZ/71 \ZZ$. We can number the rational points 
$P_i$ ($i=1,\ldots,7$) such that their images $\phi(P_i)=$ class of $(P_i-P_1)$ 
under the Abel-Jacobi map are as in the following table.

\smallskip
\vbox{
\bigskip\centerline{\def\quad{\hskip 0.6em\relax}
\def\quod{\hskip 0.5em\relax }
\vbox{\offinterlineskip
\hrule
\halign{&\vrule#&\strut\quod\hfil#\quad\cr
height2pt&\omit&&\omit&&\omit&&\omit&&\omit&&\omit&&\omit&&\omit&\cr
&$i$ && $1$ && $2$ && $3$ && $4$ && $5$ && $6$ && $7$& \cr
\noalign{\hrule}
&$\phi(P_i)$ && $0$ && $1$ && $34$ && $55$ && $10$ && $14$&& $49$& \cr
} \hrule}
}}

There is a divisor $Q_3$ of degree $3$ with $\phi(Q_3-3\, P_0)=9$
and one of degree $5$, say $Q_5$, with $\phi(Q_5-5\, P_0)=35 \, (\bmod \, 71)$;
there is also a divisor $Q_7$ of degree $7$ with $\phi(Q_7-7 P_0)=21$.
We list a number of relations:

\vbox{
\bigskip\centerline{\def\quad{\hskip 0.6em\relax}
\def\quod{\hskip 0.5em\relax }
\vbox{\offinterlineskip
\hrule
\halign{&\vrule#&\strut\quod\hfil#\quad\cr
height2pt&\omit&&\omit&&\omit&&\omit&&\omit&\cr
& n && divisor  && n && divisor &\cr
\noalign{\hrule}
& $1$ && $[-3,2,1,1,2,1,-4]$  && $5$ && $[-1, 3, 1, 1, 1, 2, -7]$  & \cr
& $2$ && $[-11,1,2,1,2,3,2]$  & & $6$ && $[2,-1,1,2,2,2,1]\sim 3Q_3$  &\cr
& $3$ && $[-13,2,1,2,2,3,3]$  && $7$ && $[1,1,1,1,1,1,1]\sim Q_7$  & \cr
& $4$ && $[1,-5,1,2,2,1,1]\sim Q_3$  && $8$ && $[-15, 4, 1, 3, 4, 1, 2]$& \cr
} \hrule}
}}

\noindent
that yield the following results:

\smallskip
\vbox{
\bigskip\centerline{\def\quad{\hskip 0.6em\relax}
\def\quod{\hskip 0.5em\relax }
\vbox{\offinterlineskip
\hrule
\halign{&\vrule#&\strut\quod\hfil#\quad\cr
height2pt&\omit&&\omit&&\omit&&\omit&&\omit&&\omit&&\omit&&\omit&\cr
& curve && $g$ && $\# C({\FF}_2)$ && interval && curve && $g$ && $\# C({\FF}_2)$ && interval & \cr
\noalign{\hrule}
& $\langle f_1\rangle $ && $9$ && $12$ && $[12]$ && $\langle f_2, f_3\rangle $ && $29$ && $25$ && $[25-27]$  &\cr
& $\langle f_7\rangle $ && $12$ && $14$ && $[14-15]$ && $\langle f_4, f_6\rangle $ && $31$ && $25$ && $[27-28]$  &\cr
& $\langle f_1, f_5 \rangle $ && $24$ && $22$ && $[21-23]$   
&& $\langle f_2,f_3,f_8\rangle$ && $69$ && $49$ && $[49-52]$ &\cr
} \hrule}
}}

Note that for $f_1$ we have to do Artin-Schreier reduction. If
the conductor would be $3P_0$ (resp.\ $3P_0+P_6$) then $C_{f_1}$
would give $g=7,N=11$ (resp.\ $g=8,12$), and these are impossible.
So the conductor is $3P_0+3P_6$ giving $(g=9,N=12)$, an optimal curve.
Again we find one improvement (for $g=24$) for our tables. 
The interval for genus $69$ comes from page 121 of \cite{NX5}.
\end{subsection}

\begin{subsection}{Base of Genus $5$ over ${\FF}_2$}
Consider the fibre product $C$ of $C_1$ given by $y^2+y=x^3+x$ and $C_2$ given
by $y^2+y=x^5+x^3$ over the $x$-line. 
This is an optimal curve of genus $5$ with $9$ points.
Its Jacobian is isogenous with the product of
the three Jacobians of the curves $C_1$ and $C_2$
and the curve $C_2'$ defined by $y^2+y=x^5+x$. The corresponding $L$-polynomials
are $2t^2+2t+1$, $4t^4+4t^3+2t^2+2t+1$ and $4t^4+4t^3+4t^2+2t+1$. The class
groups are $\ZZ/5\ZZ$, $\ZZ/13\ZZ$ and $\ZZ/15\ZZ$. 
The nine rational points are 
$P_{\infty}$ and eight points that can be identified by their 
$(x,y_1,y_2,y_3)$-coordinates. 
Writing the class group of $C$ as $\ZZ/65\ZZ \times \ZZ/15\ZZ$
we have the following table:

\begin{footnotesize}
\smallskip
\vbox{
\bigskip\centerline{\def\quad{\hskip 0.6em\relax}
\def\quod{\hskip 0.5em\relax }
\vbox{\offinterlineskip
\hrule
\halign{&\vrule#&\strut\quod\hfil#\quad\cr
height2pt&\omit&&\omit&&\omit&&\omit&&\omit&&\omit&&
\omit&&\omit&&\omit&&\omit&\cr
&$i$ && $1$ && $2$ && $3$ && $4$ && $5$ && $6$ && $7$&& $8$ && $9$ & \cr
\noalign{\hrule}
&$\phi(P_i)$&&$(0,0)$&&$(1,1)$&&$(51,14)$&& $(64,1)$&& 
$(14,14)$&&$(57,11)$&&$(47,4)$&&$(8,11)$&& $(18,4)$&\cr
} \hrule}
}}
\end{footnotesize}

The relation $[1, 2, -13, 1, 2, 2, 2, 2, 1]$ gives a curve with $g=16$
and $17$ points (interval $[17-18]$). The relation
$[1, 1, -11, 1, -1, 3, 2, 1, 3]$ gives a curve with $g=16$ and $16$
points. The fibre product has genus $g=39$ with $31$ points.
The combination of the relations $[1, 2, -13, 1, 2, 2, 2, 2, 1]$ 
and $[3,3,-17,1,2,1,3,1,3]$ gives a curve with $g=42$ and $33$ points
(interval $[33-35]$). 

\end{subsection}
\begin{subsection}{Genus $1$ over ${\FF}_3$}
Consider the elliptic curve defined by the equation $y^2=x^3+2x+1$ over
${\FF}_3$. It has $7$ points $P_i$ with $i=0,\ldots,7$ such that
$P_i \mapsto i (\bmod 7)$ defines an isomorphism $C({\FF}_3)\cong {\ZZ/7\ZZ}$.
The relations $[-4, -1, 0, 1, 2, 1, 1]$ and $[-4, -2, 1, 2, 1, 1, 1]$
give a fibre product of genus $g=30$ with $38$ points that improves 
the interval $[37-46]$ slightly.

\end{subsection}

\begin{subsection}{Curves of Genus $2$ over ${\FF}_3$}
\bigskip
The curve $C$ of genus $2$ over ${\FF}_3$ given by $y^2=x^5+x^3+x+1$
has zeta function $(9t^4+9t^3+5t^2+3t+1)/(1-t)(1-3t)$. The curve has
$7$ rational points and class group ${\ZZ}/27\ZZ$. The $7$ rational
points map to $0,1,26,17,10,23,4$ in the group. We have the relations
$[1,1,1,1,-2,2,-4]$ and $[1,1,1,2,-1,1,-5]$ and the corresponding
fibre product has $g=44$ with $47$ points. 
This gives a new entry for the table, albeit not a very strong one
(resulting interval $[47-61]$).
\end{subsection}
\bigskip
\begin{subsection}{Base Curve of Genus $3$ over ${\FF}_3$}
The curve $C$ of genus $3$ given by
$$
2\, x^4+x^3z+2\, x^2y^2+x^2yz+x^2z^2+2\, xz^3+2\, y^4+2\, y^3z+2y^2z^2=0
$$
has $10$ rational points and class group ${\ZZ}/204$. The $10$ points map
to 
$$
0,72,129,59,182,121, 172, 45, 47, 26
$$ in ${\ZZ}/204$. We have the relation
$[1, 2, 1, 1, 1, 1, 1, 1, -10, 1]$ leading to an Artin-Schreier cover with
genus $18$ with $28$ points improving the interval $[26-31]$ to $[28-31]$.
\end{subsection}
\end{section}

\begin{section}{Subgroups of the Class group}
Let $C$ be a smooth complete irreducible curve over $\FF_q$ with Jacobian
${\rm Jac}(C)$ and class group $G={\rm Jac}(C)(\FF_q)$. We choose a rational 
point $P$, provided there is one, and consider the morphism $\phi: C \to 
{\rm Jac}(C)$ given by $Q \mapsto Q-\deg(Q) P$.
If $H$ is a subgroup of $G$ of index $d$ containing the images 
$\{\phi(P_i): i \in I\}$  then there exists an unramified degree 
$d$ cover $\tilde{C}$ of $C$ in which the points $P_i$ with $i\in I$ split completely.
Choosing $C$ and $I$ appropriately can produce curves with many points.

\begin{subsection}{Base of Genus $4$ over ${\FF}_2$}
The reader might think that it is necessary to start with a
curve with many points. To dispel this idea start with the
genus $4$ curve given by the equation
$$
y^2+y=(x^7+x^5+1)/(x^2+x)
$$
over ${\FF}_2$ 
which has three rational points $P_1,P_2,P_3$. The $L$-polynomial is
$16 \, t^8+4\, t^6+4\, t^5+4\, t^4+2\, t^3+ t^2+1$, hence
the class number is $32$. In fact, the class group is
${\ZZ}/2{\ZZ} \times {\ZZ}/16{\ZZ}$ and the differences $P_i-P_1$
 map to $(0,0)$, $(1,0)$ and $(1,8)$. Hence there exists an \'etale
cover of degree $8$ in which the three points split completely.
This gives a curve of genus $25$ with $24$ points. This is optimal.

Or start with the hyperelliptic curve of genus $5$ given by
$$
y^2+y=(x^9+x^7+x^3+x+1)/(x^2+x)
$$
again with three rational points and class group ${\ZZ}/2{\ZZ} \times
{\ZZ}/30{\ZZ}$ with the points now mapping to $(0,0)$, $(1,0)$ and $(0,15)$.
This gives an \'etale cover of degree $15$ with genus $g=61$ with $45$
points, very close to the best upper bound $47$ and improving the interval
$[41-47]$ of \cite{NX5}, p.\ 121.
\end{subsection}

\begin{subsection}{Genus $2$ over ${\FF}_3$}
We consider the curve $C$ of genus $2$ over ${\FF}_3$ given by the affine
equation $y^3-y=x-1/x$. The zeta function is
$(9t^4+12t^3+10t^2+4t+1)/(1-t)(1-3t)$ and its class group is isomorphic
to ${\ZZ/6\ZZ} \times {\ZZ/6 \ZZ}$. The images of the $8$ rational
points under the Abel-Jacobi map are
$$
[0,0], [2,4], [1,0],[1,4],[0,1],[2,3],[5,5],[3,5]
$$
in $({\ZZ/6\ZZ})^2$.
The subgroup of the class group of index $12$ given by
$x_1+x_2\equiv  0 \, (\bmod \, 3)$ and $x_1\equiv x_2 \equiv  0 \,
(\bmod \, 2)$ contains $[0,0]$ and
$[2,4]$ and thus leads to a cover of degree $12$ of genus $g=13$
and with $24$ rational points (interval $[24-25]$).
There is a place $Q_2$ of degree $2$ which maps to $[2,1]$
in the class group.
The relation $\sum_{i=1}^8 P_i \sim 4Q_2$
leads to an Artin-Schreier cover with $g=14$ with $24$ points (interval
$[24-25]$).

An alternative case: consider $y^2=x(x^2+1)(x^2-x-1)$ over ${\FF}_3$ of genus
$2$ with $6$ points and class group ${\ZZ/2\ZZ}\times {\ZZ/10\ZZ}$.
The points map under an Abel-Jacobi map to
$$
[0,0],[1,0],[0,9],[0,1],[1,7],[1,3].
$$
So the index $10$ subgroup defined by $x_2=0 \, (\bmod \, 10)$ 
gives a cover of degree
$10$ of genus $11$ with $20$ points. (Interval $[20-22]$.)

Finally, consider the curve $C$ defined over ${\FF}_3$ by the
equation $y^2=2x^5+x^4+x$. It has $6$ rational points and zeta function
$(9t^4+6t^3+4t^2+2t+1)/(1-t)(1-3t)$ and class group $\ZZ/22\ZZ$. The six points
map to $0,11,12,10,3,19 \, (\bmod \, 22)$.
We thus see that there is an \'etale cover
of degree $11$ with genus $12$ and $22$ rational points (interval
$[22-23]$).

\smallskip
Let $C$ be the curve of genus $2$ over $\FF_3$ defined by the equation
$y^3-y=x+1/x$. The zeta function of this curve is
$$
(81t^4+90t^3+43t^2+10t+1)/(1-t)(1-9t)
$$
and its class group is $\ZZ/ 15 \ZZ \times \ZZ/15\ZZ$. The curve possesses $20$
rational points over ${\FF}_9$. The images of these points under a suitable Abel-Jacobi map
are given by
$$
\begin{aligned}
&[0,13],[0,8],[7,9],[8,12],[14,13],[1,8],[6,12],[9,9],[14,2],[1,4], \cr
&[2,9],[13,12],[4,7],[11,14],[9,6],[6,0],[4,8],[11,13],[0,6],[0,0].
\end{aligned}
$$
The equation $5a+5b=0 \, (\bmod \, 15)$ defines a subgroup of index $3$ in
$\ZZ/ 15 \ZZ \times \ZZ/15\ZZ$ containing the images of $10$ points,
hence we get a curve $\tilde{C}$,
an unramified cover of $C$ of degree $3$ with genus $4$ with $30$ points,
the maximum possible.

Similarly, the six points $P_i$ corresponding to pairs $[a,b]$
with $a\equiv b \equiv 0 \, (\bmod \, 3)$ lie in a subgroup of index $9$. So there is
an unramified cover $\tilde{C}$ of degree $9$ and genus $10$ with
$9\times 6=54$ rational points. The interval in the tables is $[54-55]$.

\end{subsection}
\begin{subsection}{Base Curve of Genus $3$ over ${\FF}_3$}

Consider the plane curve of degree $4$ over $\FF_3$ 
given by
$$
2\, x^3y+2\, x^3z+x^2y^2+xz^3+2\, y^3z+yz^3=0.
$$
This curve has $10$ points over $\FF_3$ and has class group isomorphic to
$({\ZZ/14\ZZ})^2$
and the $10$ points go to
$$
[0,0],[1,0],[6,13],[7,3],[13,7],[4,12],[11,11],[6,2],[4,7]]
$$
The subgroup defined by $x_2\equiv 0\, (\bmod \, 7)$ has index $7$
and contains the images of $4$ points. The corresponding curve over ${\FF}_3$
has genus $15$ with $28$ rational points; this is optimal.

Or consider the hyperelliptic curve given by the equation
$$
y^2+(x^3-x)y=x^7-x^2+x
$$
over ${\FF}_3$. It is of genus $3$, has five rational points mapping in
the class group ${\ZZ}/2{\ZZ} \times {\ZZ}/2 {\ZZ}\times {\ZZ}/12{\ZZ}$ to
$(0,0,0,), (0,0,6), (0,1,0), (0,0,11), (0,0,1)$. The 
index $12$ subgroup containing the first three points gives rise to a
cover of genus $g=25$ with $36$ points (the interval being $[36-40]$).

The hyperelliptic curve $y^2=x^7-x^2+x$ over ${\FF}_3$ has similarly
a class group ${\ZZ}/2{\ZZ}\times {\ZZ}/2{\ZZ} \times {\ZZ}/14{\ZZ}$ 
and an index $14$ subgroup containing three points giving rise to a curve
of genus $29$ with $42$ points (interval $[42-44]$).

The hyperelliptic curve of genus $4$ over ${\FF}_3$ given by $y^2+xy=x^9-x$
has class group ${\ZZ}/2{\ZZ} \times {\ZZ}/126{\ZZ}$ and the six rational
points map to 
$$
(0,0), (1,0), (0,125), (0,1), (0,14), (0,112)
$$
and so the index $14$ subgroup containing four of these
points gives rise to a cover of genus $43$ with $56$ points
improving the interval $[55-60]$.
\end{subsection}
\begin{subsection}{Examples over ${\FF}_4$}
Let $C$ be the hyperelliptic curve of genus $2$ over ${\FF}_4$ given by
$y^2+y=x^5+x^3+x$ with class group ${\ZZ/7\ZZ}\times {\ZZ/7\ZZ}$.
It has a cover of degree $7$ in which three points split completely
giving a curve of genus $8$ with $21$ rational points (interval $[21-24]$).

\end{subsection}
\begin{subsection}{Genus $3$ over ${\FF}_9$}
Let ${\alpha}$ be a generator of the multiplicative group ${\FF}_9^*$
and consider the curve $C$ of genus $3$ given by $y^3-y=\alpha^2x^4$
over the field ${\FF}_9$. This curve has $28$ rational points, the
bitangent points of a plane model. The class group is of the form
$({\ZZ}/4{\ZZ})^6$. The $28$ points map under an Abel-Jacobi map
to the following elements in $({\ZZ}/4{\ZZ})^6$.
\begin{footnotesize}
$$
\begin{matrix}
[0,0,0,0,0,0], [3,1,3,3,0,1],[3,2,1,0,3,0],[2,0,3,1,1,3],[2,1,1,1,1,3], \cr
[2,2,0,1,1,3],[3,1,0,1,1,3],[1,1,0,0,1,3],[2,1,0,2,1,3],[1,3,3,2,0,3], \cr
[1,2,2,3,2,0],[0,2,3,2,1,2],[2,0,0,0,1,2],[1,1,3,1,2,0],[3,2,1,2,0,3], \cr
[0,0,1,3,2,1],[0,3,1,1,2,2],[2,0,2,3,3,2],[3,3,0,3,2,3],[3,1,2,0,1,0], \cr
[0,3,2,0,0,2],[1,1,1,2,1,1],[1,2,0,3,3,3],[0,0,3,2,3,1],[2,1,0,1,1,0], \cr
[2,1,0,1,0,2],[2,1,0,1,2,3],[2,1,0,1,1,3]
\end{matrix}
$$
\end{footnotesize}
The index $2$ subgroup of the class group defined by the equation
$2x_1 \equiv 0 \, (\bmod \, 4)$ contains the images of $16$ points.
The corresponding \'etale covering of $C$ of degree $2$ is of genus 
$5$ and has $32$ points (interval $[32-35]$).

The index $4$ subgroup of the class group defined by the equation
$x_2+x_3+x_4+x_5 \equiv 0 \, (\bmod \, 4)$ contains the images of $12$ points.
We thus find a degree $4$ \'etale
cover of $C$ of genus $9$ with $48$ points (interval $[48-50]$).  
\end{subsection}
\end{section}
\begin{section}{Using the class group over extension fields}
Let $C$ a smooth projective curve defined over $\FF_q$ of genus $g$
and with $m$ rational points. Let $J$ be its 
Jacobian variety. It is easy to see that if $Z(C,t)$ is the zeta function
of $C/\FF_q$ and $Z_n(C,t)$ the zeta function of $C$ considered
as a curve over $\FF_{q^n}$ then we have 
$Z_n(C,t^n)=\prod_{\zeta} Z(C,\zeta t)$,
where the product is taken over the $n$th roots $\zeta$ of $1$.
If we write $Z(C,t)=L(C,t)/(1-t)(1-qt)$  and $Z_n(C,t)=L_n(C,t)/(1-t)(1-q^nt)$
we get 
$L_n(C,1)=\prod_{\zeta} L(C,\zeta)$. Moreover, we know that $\#J(\FF_q)=L(C,1)$
and $\# J(\FF_{q^n})=L_n(C,1)$ and we thus get 
$$
[J(\FF_{q^n}): J(\FF_q)]= \prod_{\zeta^n=1, \zeta \neq 1} L(C,\zeta).
$$
Since under the Abel-Jacobi map the $\FF_q$-rational points of $C$ 
map to $J(\FF_q)$ we conclude that 
there exists an unramified cover of $C$ defined over $\FF_{q^n}$
of degree $d= \prod_{\zeta^n=1, \zeta \neq 1} L(C,\zeta)$ in which all
the $\FF_q$-rational points split completely. Thus we find a curve $\tilde{C}$
of genus $d(g-1)+1$ with at least $dm$ rational points. This idea was
exploited very succesfully by Niederreiter and Xing in their papers
\cite{NX1} up to \cite{NX5}. Here we employ it  more systematically
by going through all isomorphism classes of curves of low genera.
This improves the tables at many places. 

We now list all possible $L$-functions of genus $2$ curves over $\FF_2$.
By $N_i$ we mean $\# C({\FF}_{2^i})$.

\begin{footnotesize}

\smallskip
\vbox{
\bigskip\centerline{\def\quad{\hskip 0.6em\relax}
\def\quod{\hskip 0.5em\relax }
\vbox{\offinterlineskip
\hrule
\halign{&\vrule#&\strut\quod\hfil#\quad\cr
height2pt&\omit&&\omit&&\omit&&\omit&\cr
&$f$ && $[N_1,N_2]$ && $Cl$ && $L$  &\cr
\noalign{\hrule}
height2pt&\omit&&\omit&&\omit&&\omit&\cr
&$(x^2+x)/(x^3+x^2+1)$ && $[6,6]$ && $\ZZ/19\ZZ$ && $4t^4+6t^3+5t^2+3t+1$ & \cr
&$(x^3+x+1)/(x^3+x^2+1)$ && $[0,6]$ && $(0)$ && $4t^4-6t^3+5t^2-3t+1$ & \cr
&$1/(x^3+x^2+1)$ && $[2,6]$ && $\ZZ/3\ZZ$ && $4t^4-2t^3+t^2-t+1$ &\cr
&$x/(x^3+x^2+1)$ && $[4,6]$ && $\ZZ/9\ZZ$ && $4t^4+2t^3+t^2+t+1$ & \cr
&$x^2/(x^3+x^2+1)$ && $[4,10]$ && $\ZZ/11\ZZ$ && $4t^4+2t^3+3t^2+t+1$ & \cr
&$(x^3+1)/(x^3+x^2+1)$ && $[2,10]$ && $\ZZ/5\ZZ$ && $4t^4-2t^3+3t^2+t+1$ & \cr
&$1/x(x^2+x+1)$ && $[3,7]$ && $\ZZ/6\ZZ$ && $4t^4+t^2+1$ & \cr
&$(x+1)/x(x^2+x+1)$ && $[5,7]$ && $\ZZ/14\ZZ$ && $4t^4+4t^3+3t^2+2t+1$ & \cr
&$(x^3+x^2+1)/(x(x^2+x+1)$ && $[1,7]$ && $\ZZ/2\ZZ$ && $4t^4-4t^3+3t^2-2t+1$ & \cr
&$(x^3+x^2+1)/x(x+1)$ && $[3,3]$ && $\ZZ/2\ZZ \times \ZZ/2 \ZZ$ && $4t^4-t^2+1$ & \cr
&$1/x+x^3$ && $[4,4]$ && $\ZZ/8\ZZ$ && $4t^4+2t^3+t+1$ & \cr
&$1/x+x^2+x^3$ && $[2,8]$ && $\ZZ/4\ZZ$ && $4t^4-2t^3+2t^2-t+1$ & \cr
&$1/x+1+x^3$ && $[2,4]$ && $\ZZ/2\ZZ$ && $t^4-2t^3-t+1$ & \cr
&$1/x+1+x^2+x^3$ && $[4,8]$ && $\ZZ/10\ZZ$ && $4t^4+2t^3+2t^2+t+1$ & \cr
&$x^5$ && $[3,5]$ && $\ZZ/5\ZZ$ && $4t^4+1$ & \cr
&$x^5+x^3+x$ && $[3,9]$ && $\ZZ/7\ZZ$ && $4t^4+2t^2+1$ & \cr
&$x^5+x$ && $[5,9]$ && $\ZZ/15\ZZ$ && $4t^4+4t^3+4t^2+2t+1$ & \cr
&$x^5+x+1$ && $[1,5]$ && $\ZZ/3\ZZ$ && $4t^4-4t^3+4t^2-2t+1$ & \cr
&$x^5+x^3$ && $[5,5]$ && $\ZZ/13\ZZ$ && $4t^4+4t^3+2t^2+2t+1$ & \cr
&$x^5+x^3+1$ && $[1,9]$ && $(0)$ && $4t^4-4t^3+2t^2-2t+1$ & \cr
} \hrule}
}}

\end{footnotesize}
If we now take $n=2$ we get curves defined over ${\FF}_4$; for example
we find the following cases that do not give improvements of the table
for ${\FF}_4$ but good realizations of the present records.

\begin{footnotesize}
\smallskip
\vbox{
\bigskip\centerline{\def\quad{\hskip 0.6em\relax}
\def\quod{\hskip 0.5em\relax }
\vbox{\offinterlineskip
\hrule
\halign{&\vrule#&\strut\quod\hfil#\quad\cr
height2pt&\omit&&\omit&&\omit&&\omit&\cr
&$[N_1,N_2]$ && $g$ && $\#C(\FF_4)$&& interval &\cr
\noalign{\hrule}
&$[5,9]$ && $4$ && $15$ && $[15]$&\cr
&$[4,10]$ && $6$ && $20$ && $[20]$&\cr
&$[3,9]$ && $8$ && $21$ && $[21-24]$&\cr
} \hrule}
}}
\end{footnotesize}

With $n=3$ we get the following cases

\begin{footnotesize}
\smallskip
\vbox{
\bigskip\centerline{\def\quad{\hskip 0.6em\relax}
\def\quod{\hskip 0.5em\relax }
\vbox{\offinterlineskip
\hrule
\halign{&\vrule#&\strut\quod\hfil#\quad\cr
height2pt&\omit&&\omit&&\omit&&\omit&\cr
&$[N_1,N_2]$ && $g$ && $\#C(\FF_8)$&& interval &\cr
\noalign{\hrule}
&$[5,7]$ && $8$ && $35$ && $[35-42]$&\cr
&$[5,5]$ && $14$ && $65$ && $[65]$&\cr
&$[4,4]$ && $20$ && $76$ && $[68-83]$&\cr
} \hrule}
}}
\end{footnotesize}

Here the case $g=14$ with $\#C({\FF}_4)=65$ is optimal; the case $g=8$ with 
$35$ rational points
equals the lower bound of the interval $[35-42]$ of the tables while
$g=20$ with $76$ points improves the interval $[68-83]$ of the tables.

Starting with the curves of genus $2$ over $\FF_4$ gives the following results:

\begin{footnotesize}
\smallskip
\vbox{
\bigskip\centerline{\def\quad{\hskip 0.6em\relax}
\def\quod{\hskip 0.5em\relax }
\vbox{\offinterlineskip
\hrule
\halign{&\vrule#&\strut\quod\hfil#\quad\cr
height2pt&\omit&&\omit&&\omit&&\omit&\cr
&$[N_1,N_2]$ && $g$ && $\#C(\FF_{16})$ && interval &\cr
\noalign{\hrule}
&$[9,24]$ && $9$ && $72$ && $[72-81]$&\cr
&$[9,25]$ && $10$ && $81$ && $[81-87]$&\cr
&$[8,24]$ && $11$ && $80$ && $[80-91]$&\cr
&$[8,26]$ && $12$ && $88$ && $[83-97]$&\cr
&$[7,27]$ && $15$ && $98$ && $[98-113]$&\cr
&$[7,31]$ && $17$ && $112$ && $[112-123]$&\cr
} \hrule}
}}
\end{footnotesize}

And the same can be done starting with the curves of genus $2$ over $\FF_8$.
Here there are still many empty places in the tables due to the lack of
good examples.

\begin{footnotesize}
\smallskip
\vbox{
\bigskip\centerline{\def\quad{\hskip 0.6em\relax}
\def\quod{\hskip 0.5em\relax }
\vbox{\offinterlineskip
\hrule
\halign{&\vrule#&\strut\quod\hfil#\quad\cr
height2pt&\omit&&\omit&&\omit&&\omit&\cr
&$[N_1,N_2]$ && $g$ && $\#C(\FF_{64})$ && interval &\cr
\noalign{\hrule}
&$[18,54]$ && $20$ && $342$ && &\cr
&$[17,63]$ && $25$ && $408$ && &\cr
&$[17,65]$ && $26$ && $425$ && &\cr
&$[16,64]$ && $27$ && $416$ && &\cr
&$[16,70]$ && $30$ && $464$ && &\cr
&$[15,71]$ && $33$ && $480$ && &\cr
&$[14,66]$ && $34$ && $462$ && &\cr
&$[15,75]$ && $35$ && $510$ && &\cr
&$[14,70]$ && $36$ && $490$ && &\cr
&$[15,79]$ && $37$ && $540$ && &\cr
&$[14,74]$ && $38$ && $518$ && &\cr
&$[13,67]$ && $39$ && $494$ && $[489-650]$ &\cr
&$[14,78]$ && $40$ && $546$ && &\cr
&$[14,80]$ && $41$ && $560$ && &\cr
&$[14,82]$ && $42$ && $574$ && &\cr
&$[12,66]$ && $44$ && $516$ && &\cr
&$[13,79]$ && $45$ && $572$ && &\cr
&$[13,81]$ && $46$ && $585$ && &\cr
&$[13,83]$ && $47$ && $598$ && &\cr
&$[12,74]$ && $48$ && $564$ && &\cr
&$[13,87]$ && $49$ && $624$ && &\cr
&$[12,78]$ && $50$ && $588$ && &\cr
} \hrule}
}}
\end{footnotesize}

We can use curves of higher genus too. 
Going through all non-hyperelliptic curves of genus $3$ 
over ${\FF}_2$
gives two improvements for the table over ${\FF}_4$:

\begin{footnotesize}
\smallskip
\vbox{
\bigskip\centerline{\def\quad{\hskip 0.6em\relax}
\def\quod{\hskip 0.5em\relax }
\vbox{\offinterlineskip
\hrule
\halign{&\vrule#&\strut\quod\hfil#\quad\cr
height2pt&\omit&&\omit&&\omit&&\omit&\cr
&$[N_1,N_2,N_3]$ && $g$ && $\#C(\FF_4)$ && interval &\cr
\noalign{\hrule}
&$[5,9,5]$ && $15$ && $35$ && $[33-37]$&\cr
&$[5,11,5]$ && $17$ && $40$ && $[40]$&\cr
&$[4,12,7]$ && $27$ && $52$ && $[50-56]$&\cr
} \hrule}
}}
\end{footnotesize}
\noindent
and considering these curves over ${\FF}_4$ gives
three improvements for the table over ${\FF}_{16}$:

\begin{footnotesize}
\smallskip
\vbox{
\bigskip\centerline{\def\quad{\hskip 0.6em\relax}
\def\quod{\hskip 0.5em\relax }
\vbox{\offinterlineskip
\hrule
\halign{&\vrule#&\strut\quod\hfil#\quad\cr
height2pt&\omit&&\omit&&\omit&&\omit&\cr
&$[N_1,N_2,N_3]$ && $g$ && $\#C(\FF_{16})$ && interval &\cr
\noalign{\hrule}
&$[4,12,7]$ && $27$ && $156$ && $[145-176]$&\cr
&$[3,11,6]$ && $35$ && $187$ && &\cr
&$[5,11,5]$ && $41$ && $220$ && $[216-249]$&\cr
} \hrule}
}}
\end{footnotesize}

Going through the hyperelliptic curves of genus $4$ over ${\FF}_2$
yields one improvement of the table over ${\FF}_4$: the curve
with $[N_1,N_2,N_3,N_4]=[6,10,6,26]$ yields a curve of genus $28$ with $54$
points over ${\FF}_4$; the interval was $[53-58]$.

Taking the curve of genus
$4$ over ${\FF}_2$ with $[N_1,N_2,N_3,N_4]=[8,8,8,16]$ over ${\FF}_2$ and
$L$-polynomial $16t^8+40t^7+56t^6+56t^5+44t^4+28t^3+14t^2+5t+1$
gives a degree $13$ cover with genus $40$ and $104$ points,
improving the interval $[103-141]$.

Similarly, we can list all pairs $[N_1,N_2]$ with $N_i=\# C({\FF}_{3^i})$
occuring for curves of genus $2$ over $\FF_3$. This yields the following
harvest. 

\begin{footnotesize}
\smallskip
\vbox{
\bigskip\centerline{\def\quad{\hskip 0.6em\relax}
\def\quod{\hskip 0.5em\relax }
\vbox{\offinterlineskip
\hrule
\halign{&\vrule#&\strut\quod\hfil#\quad\cr
height2pt&\omit&&\omit&&\omit&&\omit&\cr
&$[N_1,N_2]$ && $g$ && $\#C(\FF_9)$ && interval &\cr
\noalign{\hrule}
&$[8,14]$ && $5$ && $32$ && $[32-35]$&\cr
&$[7,15]$ && $6$ && $35$ && $[35-40]$&\cr
&$[6,16]$ && $8$ && $42$ && $[40-47]$&\cr
&$[6,18]$ && $9$ && $48$ && $[48-50]$&\cr
&$[6,20]$ && $10$ && $54$ && $[54]$&\cr
} \hrule}
}}
\end{footnotesize}

Doing the same for curves of genus $2$ over ${\FF}_9$ gives the following.

\begin{footnotesize}
\smallskip
\vbox{
\bigskip\centerline{\def\quad{\hskip 0.6em\relax}
\def\quod{\hskip 0.5em\relax }
\vbox{\offinterlineskip
\hrule
\halign{&\vrule#&\strut\quod\hfil#\quad\cr
height2pt&\omit&&\omit&&\omit&&\omit&\cr
&$[N_1,N_2]$ && $g$ && $\#C(\FF_{81})$ && interval &\cr
\noalign{\hrule}
&$[20-68]$ && $26$ && $500$ && &\cr
&$[19-75]$ && $30$ && $551$ && &\cr
&$[18-78]$ && $33$ && $576$ && &\cr
&$[18-80]$ && $34$ && $594$ && $[494-689]$ &\cr
&$[18-82]$ && $35$ && $612$ && &\cr
&$[18-86]$ && $37$ && $648$ && $[568-742]$ &\cr
&$[17-83]$ && $38$ && $629$ && &\cr
&$[17-85]$ && $39$ && $646$ && &\cr
&$[17-87]$ && $40$ && $663$ && &\cr
&$[17-89]$ && $41$ && $680$ && &\cr
&$[17-91]$ && $42$ && $697$ && &\cr
&$[16-86]$ && $43$ && $672$ && &\cr
&$[16-90]$ && $45$ && $704$ && &\cr
&$[16-92]$ && $46$ && $720$ && &\cr
&$[15-85]$ && $47$ && $690$ && &\cr
&$[16-96]$ && $48$ && $752$ && $[676-885]$ &\cr
&$[16-86]$ && $49$ && $672$ && $[656-898]$ &\cr
&$[16-100]$ && $50$ && $784$ && &\cr
} \hrule}
}}
\end{footnotesize}

Similarly, the genus $2$ curves over $\FF_3$ give with $n=3$:

\begin{footnotesize}
\smallskip
\vbox{
\bigskip\centerline{\def\quad{\hskip 0.6em\relax}
\def\quod{\hskip 0.5em\relax }
\vbox{\offinterlineskip
\hrule
\halign{&\vrule#&\strut\quod\hfil#\quad\cr
height2pt&\omit&&\omit&&\omit&&\omit&\cr
&$[N_1,N_2]$ && $g$ && $\#C(\FF_{27})$ && interval & \cr
\noalign{\hrule}
&$[8,10]$ && $26$ && $200$ &&&\cr
&$[7,13]$ && $29$ && $196$ &&&\cr
&$[7,11]$ && $40$ && $273$ && $[244-346]$ &\cr
} \hrule}
}}

\end{footnotesize}

\bigskip

For the case of non-hyperelliptic curves of genus $3$ over $\FF_3$
we find the following cases: 

\begin{footnotesize}
\smallskip
\vbox{
\bigskip\centerline{\def\quad{\hskip 0.6em\relax}
\def\quod{\hskip 0.5em\relax }
\vbox{\offinterlineskip
\hrule
\halign{&\vrule#&\strut\quod\hfil#\quad\cr
height2pt&\omit&&\omit&&\omit&&\omit&&\omit&&\omit&\cr
&$m$ && $\# C(\FF_3)$ && $L(-1)$ && $g$ && $\# \tilde{C}(\FF_9)$&&{\rm interval}& \cr
\noalign{\hrule}
&$687439$ && $8$   && $12$ && $25$ &&
$96$&& $[82-108]$ & \cr
&$787452$ && $8$ && $13$ && $27$ && $104$&& $[91-114]$ & \cr
&$687411$ && $7$   && $17$ && $35$ && $119$&& $[116-139]$ &\cr
&$787567$ && $7$    && $18$ && $37$ && $126$&& $[120-145]$ &\cr
&$884286$ && $7$    && $20$ && $41$ && $140$&& $[128-158]$ &\cr
&$687541$ && $7$ && $21$ && $43$ && $147$&& $[120-164]$ & \cr
} \hrule}
}}
\end{footnotesize}

Here we use the following notation for plane curves of degree $3$
over ${\FF}_3$. We write the polynomials in $x,y,z$ in pure lexicographic
order and use then the coefficients $c_i$ ($i=1,\ldots,15)$ to
associate the natural number $m=\sum_{i=1}^{15} c_i 3^{i-1}$ to the curve.
So $x^4$ corresponds to $1$ and the Klein curve
$x^3y+y^3z+z^3x=0$ to $196833$.

\end{section}
\vfill
\eject

\begin{section}{Tables}

\font\tablefont=cmr8
\def\quad{\hskip 0.6em\relax}
\def\quod{\hskip 0.6em\relax}
\def\vhop{
    height2pt&\omit&&\omit&&\omit&&\omit&&\omit&&\omit&&\omit&&\omit&\cr}
\noindent{\bf Table p=2.}
$$
\vcenter{
\tablefont
\lineskip=1pt
\baselineskip=10pt
\lineskiplimit=0pt
\setbox\strutbox=\hbox{\vrule height .7\baselineskip
                                depth .3\baselineskip width0pt}%
\offinterlineskip
\hrule
\halign{&\vrule#&\strut\quod\hfil#\quad\cr
\vhop
&$g\backslash q$&&2&&4&&8&&16&&32&&64&&128&\cr
\vhop
\noalign{\hrule}
\vhop
&1&&5&&9&&14&&25&&44&&81&&150&\cr
&2&&6&&10&&18&&33&&53&&97&&172&\cr
&3&&7&&14&&24&&38&&64&&113&&192&\cr
&4&&8&&15&&25&&45&&71--74&&129&&215&\cr
&5&&9&&17&&29--30&&49--53&&83--85&&132--145&&227--234&\cr
&6&&10&&20&&33--35&&65&&86--96&&161&&243--258&\cr
&7&&10&&21--22&&34--38&&63--69&&98--107&&177&&262--283&\cr
&8&&11&&21--24&&35--42&&62--75&&97--118&&169--193&&276--302&\cr
&9&&12&&26&&45&&72--81&&108--128&&209&&288--322&\cr
&10&&13&&27&&42--49&&81--87&&113--139&&225&&296--345&\cr
\noalign{\hrule}
&11&&14&&26--29&&48--53&&80--91&&120--150&&201--236&&294--366&\cr
&12&&14--15&&29--31&&49--57&&{\color{red}88}--97&&129--161&&257&&321--388&\cr
&13&&15&&33&&56--61&&97--102&&129--172&&225--268&&&\cr
&14&&15--16&&32--35&&65&&97--107&&146--183&&241--284&&353--437&\cr
&15&&17&&35--37&&57--67&&98--113&&158--194&&258--300&&386--455&\cr
&16&&17--18&&36--38&&56--71&&95--118&&147--204&&267--316&&&\cr
&17&&17--18&&40&&63--74&&112--123&&154--212&&&&&\cr
&18&&18--19&&41--42&&65--77&&113--129&&161--220&&281--348&&&\cr
&19&&20&&37--43&&60--80&&129--134&&172--228&&315--364&&&\cr
&20&&19--21&&40--45&&{\color{red}76}--83&&127--139&&177--236&&{\color{red}342}--380&&&\cr
\noalign{\hrule}
&21&&21&&44--47&&72--86&&129--145&&185--243&&281--396&&&\cr
&22&&21--22&&42--48&&74--89&&129--150&&&&321--412&&&\cr
&23&&22--23&&45--50&&68--92&&126--155&&&&&&&\cr
&24&&{\color{red}22}--23&&49--52&&81--95&&129--161&&225--267&&337--444&&513--653&\cr
&25&&24&&51--53&&86--97&&144--165&&&&{\color{red}408}--460&&&\cr
&26&&24--25&&55&&82--100&&150--171&&&&{\color{red}425}--476&&&\cr
&27&&24--25&&52--56&&96--103&&156--176&&213--290&&{\color{red}416}--492&&&\cr
&28&&25--26&&54--58&&97--106&&145--181&&257--298&&513&&577--745&\cr
&29&&25--27&&52--60&&97--109&&161--186&&227--305&&&&&\cr
&30&&25--27&&53--61&&96--112&&162--191&&273--313&&{\color{red}464}--535&&609--784&\cr
\noalign{\hrule}
&31&&27--28&&60--63&&89--115&&165--196&&&&450--547&&578--807&\cr
&32&&{\color{red}27}--29&&57--65&&90--118&&&&&&&&&\cr
&33&&28--29&&65--66&&97-121&&193--207&&&&{\color{red}480}--570&&&\cr
&34&&27--30&&65--68&&98--124&&183--212&&&&{\color{red}462}--581&&&\cr
&35&&29--31&&64--69&&112--127&&187--217&&253--351&&{\color{red}510}--593&&&\cr
&36&&30--31&&64--71&&112--130&&185--222&&&&{\color{red}490}--604&&705--917&\cr
&37&&30--32&&66--72&&121--132&&208--227&&&&{\color{red}540}--616&&&\cr
&38&&30--33&&64--74&&129--135&&193--233&&291--375&&{\color{red}518}--627&&&\cr
&39&&33&&65--75&&120--138&&194--238&&&&{\color{red}494}--638&&&\cr
&40&&32--34&&75--77&&103--141&&225--243&&293--390&&{\color{red}546}--649&&&\cr
\noalign{\hrule}
&41&&33--35&&65--78&&118--144&&220--249&&308--398&&{\color{red}560}--661&&&\cr
&42&&33--35&&75--80&&129--147&&209--254&&307--405&&{\color{red}574}--672&&&\cr
&43&&{\color{red}34}--36&&72--81&&116--150&&226--259&&306--413&&{\color{red}546}--684&&&\cr
&44&&33--37&&68--83&&130--153&&226--264&&325--420&&{\color{red}516}--695&&&\cr
&45&&33--37&&80--84&&144--156&&242--268&&313--428&&{\color{red}572}--706&&&\cr
&46&&34--38&&81--86&&129--158&&243--273&&&&{\color{red}585}--717&&&\cr
&47&&36--38&&73--87&&126--161&&&&&&{\color{red}598}--729&&&\cr
&48&&34--39&&80--89&&128--164&&243--282&&&&{\color{red}564}--740&&&\cr
&49&&36--40&&81--90&&130--167&&213--286&&&&{\color{red}624}--751&&913--1207&\cr
&50&&40&&91--92&&130--170&&255--291&&&&{\color{red}588}--762&&&\cr
\vhop
}
\hrule
}
$$
%
\vfill\eject
\font\tablefont=cmr8
\def\quad{\hskip 0.6em\relax}
\def\quod{\hskip 0.6em\relax}
\def\vhop{height2pt&\omit&&\omit&&\omit&&\omit&&\omit&\cr}
\noindent{\bf Table p=3.}
\bigskip
$$
\vcenter{
\tablefont
\lineskip=1pt
\baselineskip=10pt
\lineskiplimit=0pt
\setbox\strutbox=\hbox{\vrule height .7\baselineskip
                                depth .3\baselineskip width0pt}%
\offinterlineskip
\hrule
\halign{&\vrule#&\strut\quod\hfil#\quad\cr
height2pt&\omit&&\omit&&\omit&&\omit&&\omit&\cr
&$g\backslash q$&&3&&9&&27&&81&\cr
height2pt&\omit&&\omit&&\omit&&\omit&&\omit&\cr
\vhop
\noalign{\hrule}
\vhop
height2pt&\omit&&\omit&&\omit&&\omit&&\omit&\cr
&1&&7&&16&&38&&100&\cr
&2&&8&&20&&48&&118&\cr
&3&&10&&28&&56&&136&\cr
&4&&12&&30&&64&&154&\cr
&5&&13&&32--35&&72--75&&160--172&\cr
&6&&14&&35--40&&76--85&&190&\cr
&7&&16&&40--43&&82--95&&180--208&\cr
&8&&17--18&&{\color{red}42}--47&&92--105&&226&\cr
&9&&19&&48--50&&99--113&&244&\cr
&10&&20--21&&54&&94--123&&226--262&\cr
\noalign{\hrule}
&11&&20--22&&55--58&&100--133&&220--280&\cr
&12&&22--23&&56--62&&109--143&&298&\cr
&13&&24--25&&64--65&&136--153&&256--312&\cr
&14&&24--26&&56--69&&&&278--330&\cr
&15&&28&&64--73&&136--170&&292--348&\cr
&16&&27--29&&74--77&&144--178&&370&\cr
&17&&25--30&&74--81&&128--185&&288--384&\cr
&18&&{\color{red}28}--31&&67--84&&148--192&&306--401&\cr
&19&&32&&84--88&&145--199&&&\cr
&20&&30--34&&70--91&&&&&\cr
\noalign{\hrule}
\noalign{\hrule}
&21&&32--35&&88--95&&163--213&&352--455&\cr
&22&&30--36&&78--98&&&&&\cr
&23&&32--37&&92--101&&&&&\cr
&24&&31--38&&91--104&&208--234&&&\cr
&25&&36--40&&{\color{red}96}--108&&196--241&&392--527&\cr
&26&&36--41&&110--111&&{\color{red}200}--248&&{\color{red}500}--545&\cr
&27&&39--42&&{\color{red}104}--114&&&&&\cr
&28&&37--43&&105--117&&&&&\cr
&29&&42--44&&104--120&&{\color{red}196}--269&&&\cr
&30&&{\color{red}38}--46&&91--123&&196--276&&{\color{red}551}--617&\cr
\noalign{\hrule}
&31&&40--47&&120--127&&&&460--635&\cr
&32&&40--48&&92--130&&&&&\cr
&33&&46--49&&128--133&&220--297&&{\color{red}576}--671&\cr
&34&&{\color{red}46}--50&&111--136&&&&{\color{red}594}--689&\cr
&35&&47--51&&{\color{red}119}--139&&&&{\color{red}612}--707&\cr
&36&&48--52&&118--142&&244--318&&730&\cr
&37&&52--54&&{\color{red}126}--145&&236--325&&{\color{red}648}--742&\cr
&38&&&&105--149&&&&629--755&\cr
&39&&48--56&&140--152&&271--340&&{\color{red}646}--768&\cr
&40&&56--57&&118--155&&{\color{red}273}--346&&{\color{red}663}--781&\cr
\noalign{\hrule}
&41&&50--58&&{\color{red}140}--158&&&&{\color{red}680}--795&\cr
&42&&52--59&&122--161&&280--360&&{\color{red}697}--808&\cr
&43&&{\color{red}56}--60&&{\color{red}147}--164&&&&{\color{red}672}--821&\cr
&44&&{\color{red}47}--61&&119--167&&278--374&&&\cr
&45&&54--62&&136--170&&&&{\color{red}704}--847&\cr
&46&&55--63&&162--173&&&&{\color{red}720}--859&\cr
&47&&54--65&&154--177&&299--395&&{\color{red}690}--872&\cr
&48&&55--66&&163--180&&325--402&&{\color{red}752}--885&\cr
&49&&64--67&&168--183&&316--409&&768--898&\cr
&50&&63--68&&182--186&&312--416&&{\color{red}784}--911&\cr
\vhop
}
\hrule
}
$$
\eject
For exhaustive references we refer to the bibliography of the tables, 
\cite{Tables}. For general background on curves over finite fields we refer
to the book by Stichtenoth \cite{St1} and for general background on curves
over finite fields with many points to Serre's notes \cite{S2} and for
an overview of the methods of Niederreiter and Xing to \cite{NX5}.

\end{section}

\end{document}